\documentclass{article}

\usepackage{arxiv}
\usepackage[utf8]{inputenc} 
\usepackage[T1]{fontenc}    
\usepackage{hyperref}       
\usepackage{url}            
\usepackage{booktabs}       
\usepackage{amsfonts}       
\usepackage{nicefrac}       
\usepackage{microtype}      
\usepackage{amsmath}
\usepackage{graphicx}
\usepackage{doi}
\graphicspath{ {./images/} }

\hypersetup{
pdftitle={A Numerical Technique for Coupling the Momentum and the Continuity Equations for Semi-Implicit 3D Ocean Models},
pdfsubject={physics.flu-dyn, math.NA},
pdfauthor={Ali Shahabi, Reza Ghiassi},
pdfkeywords={free-surface flow, shallow water, finite difference, semi-implicit, hydrostatic pressure, three-dimensional},
}

\title{A Numerical Technique for Coupling the Momentum and the Continuity Equations for Semi-Implicit 3D Ocean Models}

\author{
 Ali Shahabi \\
  Department of Civil and Environmental Engineering\\
  Old Dominion University, Norfolk, VA\\
  \texttt{ashah010@odu.edu} \\
  \And
 Reza Ghiassi \\
  Faculty of Civil Engineering, University College of Engineering\\
  University of Tehran, Tehran, Iran\\
  \texttt{rghiassi@ut.ac.ir} \\
}

\begin{document}
\maketitle

\begin{abstract}
Semi-implicit methods are powerful and efficient tools for the three-dimensional modeling of coastal and oceanic processes. A semi-implicit finite difference method for 3D hydrostatic primitive equations is presented in this paper. The governing equations are time-discretized with an implicit treatment of barotropic pressure gradient and vertical viscosity. The discretized momentum equations along a water column are coupled with the depth-integrated continuity equation of the column to construct a linear system for free-surface elevations. The novelty of this work lies in formulating an efficient method with the order of complexity O(N) for coupling the momentum and the continuity equations. In this method, the horizontal velocity components are expressed in terms of neighbor free-surface elevations by some simple recursive formulas and then are substituted in the integrated continuity equation.
\end{abstract}

\keywords{free-surface flow \and shallow water \and finite difference \and semi-implicit \and hydrostatic pressure \and three-dimensional}

\pagestyle{plain}

\section{Introduction}
Numerical modeling of free-surface flows has been widely used to study the hydrodynamics and the transport processes of oceans, coastal waters, and rivers. Many of the geophysical flow models are based on Hydrostatic Primitive Equations (HPE). The HPEs relax the assumption of hydrostatic pressure, which reduces the vertical momentum equation to the hydrostatic pressure distribution. Large-scale phenomena in the free-surface flows are accurately described by the HPEs \cite{marshall1997}. Several three-dimensional hydrostatic models have been developed and applied successfully to various practical problems. Different numerical methods, horizontal and vertical grid systems, and time discretization approaches have been employed to develop models.

In geophysical flow modelling, a reasonable balance between accuracy, stability, and computational efficiency is generally achieved when time discretization of the barotropic pressure gradient and the vertical eddy viscosity in momentum equations are implicit. This is due to the severe stability limitations that arise from these terms \cite{vreugdenhil1994}. Different solution methods have been presented yet for the HPEs.

Some models (e.g., POM \cite{blumberg1987}, ROMS \cite{shchepetkin2005}, and FVCOM \cite{chen2003}) utilize mode-splitting technique \cite{madala1977}, in which the pressure gradient terms are not often treated implicitly. However, in order to retain computational efficiency, the free-surface elevations are determined externally from the depth-integrated equations (external mode) and the velocity components are calculated from the main three-dimensional equations (internal mode). Therefore, a large time-step can be used for the internal mode because it is no longer limited by CFL stability conditions. Nevertheless, models using this technique have errors associated with using two systems of equations. Recent efforts in this field have led to the development of more schemes to solve the three-dimensional hydrostatic systems.

A free-surface correction method developed by Chen presents a predictor-corrector method for the HPEs \cite{chen2003a}. In this method, first an intermediate free-surface is obtained with the explicit and implicit treatment of the pressure gradient and the vertical viscosity terms. Then the intermediate free-surface is corrected by solving a five-diagonal linear system for free-surface changes over the time step.

A semi-implicit finite difference method has been first presented by Casulli and Cheng which solves the three-dimensional HPEs directly with an implicit treatment of the pressure gradient and the vertical viscosity \cite{casulli1992}. The semi-implicit algorithm is accurate, stable, and straightforward. It calculates the free-surface elevations and the horizontal velocity components simultaneously. This method has been also used to develop quasi-hydrostatic \cite{casulli1998}, non-hydrostatic \cite{fringer2006}, unstructured grid \cite{casulli2000}, and isopycnal coordinate \cite{casulli1997} models. In recent years, semi-implicit methods have been used in several studies. The semi-implicit models such as TRIM-3D \cite{casulli1992}, UnTRIM \cite{casulli2000}, ECOM-si \cite{blumberg1994}, SUSTAN \cite{fringer2006}, and ELCIRC \cite{zhang2004} have enjoyed great success in simulating geophysical flows.

In the semi-implicit methods, the solution approach is to couple the horizontal momentum equations with the depth-integrated continuity equations in a way that a linear system is constructed for the free-surface. For this purpose, the horizontal velocity components should be derived in terms of neighbor water elevations. By substituting them in the depth-integrated continuity equation, a five-diagonal (or sparse in case of unstructured grid) linear system for free-surface is achieved. In the discretized momentum equation, each horizontal velocity component is related with both the neighbor water levels and the neighbor horizontal velocities (along the column). A basic way to solve such a system for the horizontal velocity is to arrange all the momentum equations along a water column in a tridiagonal linear system and inverting the coefficient matrix. Doing so is however computationally expensive and unjustified in practice.

In this paper, an algorithm is presented to construct a linear system for free-surface elevations. A three-dimensional shallow water model is developed to illustrate the details and validate the algorithm. The governing equations are discretized by finite difference method. In the proposed method, the horizontal velocity components are determined in terms of the free-surface elevations by some simple recursive formulas. The derived formulas are then substituted in the depth-integrated continuity equation to calculate the free-surface elevations.

\section{Governing Equations}
The conservation laws of mass and momentum that describe the fluid motion are known as Reynolds-Averaged Navier-Stokes equations (RANS). The three-dimensional hydrostatic primitive equations (HPE) are derived from the Navier-Stokes equations by involving the Coriolis force and applying the hydrostatic approximation. Provided that the density is constant, a conservative form of the equations is expressed as:

\begin{equation}
\begin{aligned}
    &\frac{\partial u}{\partial t} + u \frac{\partial u}{\partial x} + v \frac{\partial u}{\partial y} + w \frac{\partial u}{\partial z} &= -g \frac{\partial \xi}{\partial x} + \nu_h \left( \frac{\partial^2 u}{\partial x^2} + \frac{\partial^2 u}{\partial y^2} \right) + \frac{\partial}{\partial z} \left( \nu_z \frac{\partial u}{\partial z} \right) + f v, \\
\end{aligned}
\end{equation}
\begin{equation}
\begin{aligned}
    &\frac{\partial v}{\partial t} + u \frac{\partial v}{\partial x} + v \frac{\partial v}{\partial y} + w \frac{\partial v}{\partial z} &= -g \frac{\partial \xi}{\partial y} + \nu_h \left( \frac{\partial^2 v}{\partial x^2} + \frac{\partial^2 v}{\partial y^2} \right) + \frac{\partial}{\partial z} \left( \nu_z \frac{\partial v}{\partial z} \right) - f u,\\
\end{aligned}
\end{equation}
\begin{equation}
\begin{aligned}
    &\frac{\partial u}{\partial x} + \frac{\partial v}{\partial y} + \frac{\partial w}{\partial z} &= 0.\\
\end{aligned}
\end{equation}

where \(u\), \(v\), and \(w\) are the velocity components in the \(x\), \(y\), and \(z\) directions respectively, \(t\) is the time, \(\xi\) is the free-surface elevation measured from the still water level, \(g\) is the gravity constant, \(f\) is the Coriolis parameter, and \(\nu_h\), \(\nu_z\) are the horizontal and vertical viscosity coefficients. Because of the small effects of the horizontal viscosity terms in shallow water environments, it is justified to use a constant coefficient for the horizontal viscosity \cite{vreugdenhil1994}. The vertical eddy viscosity coefficient is calculated using the Prandtl's mixing-length theory \cite{smagorinsky1963}.

The mathematical expression for the kinematic boundary condition at the impermeable bottom is:

\begin{equation}
    u \frac{\partial h}{\partial x} + v \frac{\partial h}{\partial y} + w = 0,
\end{equation}

where \(h\) is the water depth measured from the still water level. At the moving free surface, the kinematic boundary condition is:

\begin{equation}
    \frac{\partial \xi}{\partial t} + u \frac{\partial \xi}{\partial x} + v \frac{\partial \xi}{\partial y} = w.
\end{equation}

The depth-integrated continuity equation is derived by using the kinematic boundary conditions and integrating the continuity Eq. (3) over the water depth giving:

\begin{equation}
    \frac{\partial \xi}{\partial t} + \frac{\partial}{\partial x} \int_{-h}^{\xi} u \, dz + \frac{\partial}{\partial y} \int_{-h}^{\xi} v \, dz = 0.
\end{equation}

The dynamic boundary condition at the free surface is determined by wind shear stress giving:

\begin{equation}
    \nu \left. \frac{\partial u}{\partial z} \right|_{\xi} = \tau_{wx}, \quad \nu \left. \frac{\partial v}{\partial z} \right|_{\xi} = \tau_{wy},
\end{equation}

where \(\tau_{wx}\) and \(\tau_{wy}\) are the prescribed wind stress components \cite{taylor1916}. Similarly, at the impermeable bottom boundary, the dynamic boundary condition is specified by the bed shear stress and is expressed as:

\begin{equation}
    \nu \left. \frac{\partial u}{\partial z} \right|_{-h} = \tau_{bx}, \quad \nu \left. \frac{\partial v}{\partial z} \right|_{-h} = \tau_{by},
\end{equation}

where \(\tau_b\) is the bed shear stress. The bed shear stress can be represented in the form of a quadratic friction law giving:

\begin{equation}
    \tau_{bx} = \frac{g \sqrt{u^2 + v^2}}{C^2} u, \quad \tau_{by} = \frac{g \sqrt{u^2 + v^2}}{C^2} v,
\end{equation}

where \(C\) is Chezy's roughness coefficient.

\section{Description of Proposed Algorithm}
A finite difference method is used to discretize the governing equations and the boundary conditions on a staggered grid layout in a Cartesian coordinate system. The computational domain is discretized by $n_i$, $n_j$, and $n_k$ cells in the $x$, $y$, and $z$ directions, respectively. $\Delta x$, $\Delta y$, and $\Delta z$ are the cell sizes in $x$, $y$, and $z$ directions. In the $x$ and $y$ directions, a uniform cell size is used, while in the vertical direction $\Delta z$ varies with $k$. The height of top layer cells also varies with time and horizontal location to simulate the free-surface. Each cell is denoted by $(i,j,k)$ at center and the cell faces are denoted by $(i\pm 1/2,j,k)$, $(i,j\pm 1/2,k)$ and $(i,j,k\pm 1/2)$, as shown in Figure 1. The scalar quantities are defined at the cell centers, while the velocity components are located at the faces.

\begin{figure}[ht]
\centering
\includegraphics[width=0.8\textwidth]{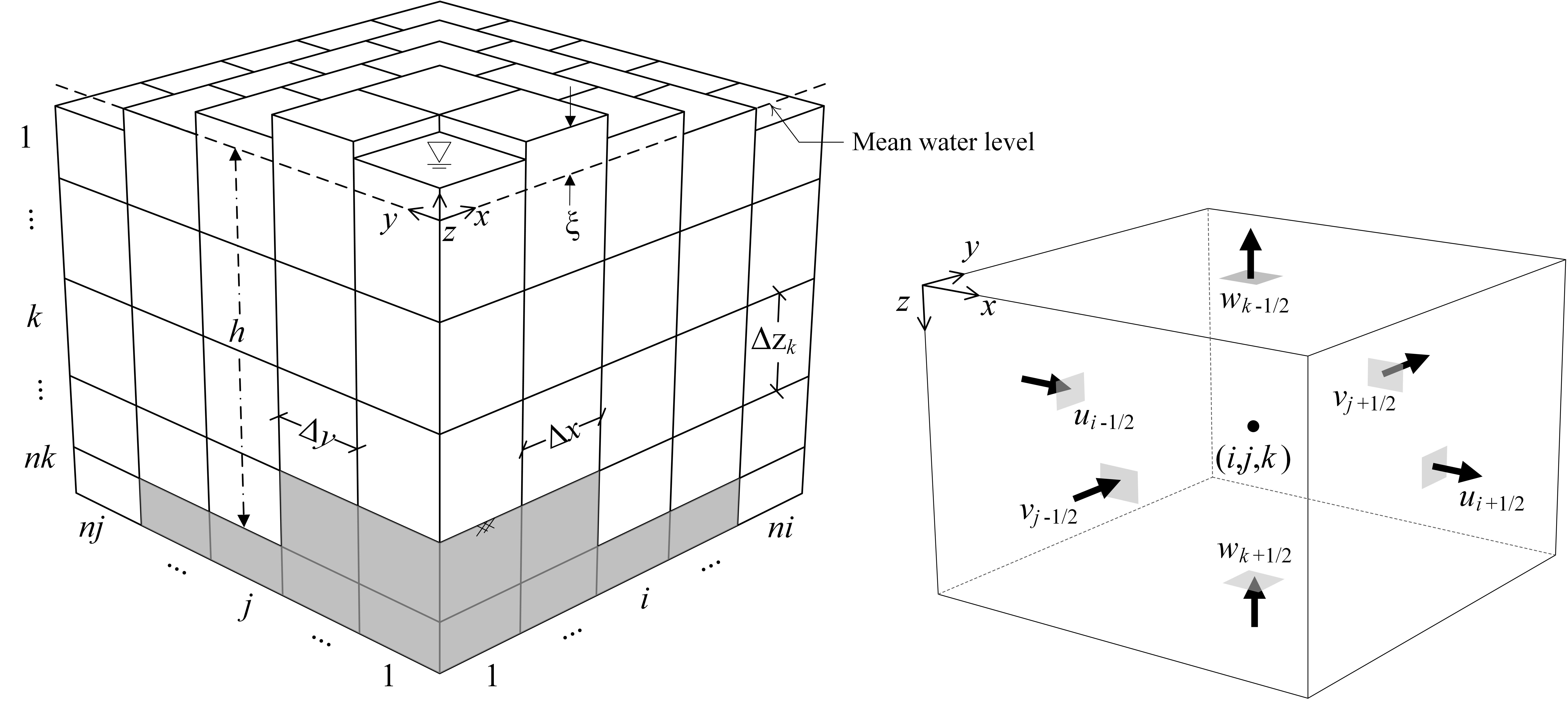}
\caption{Computational mesh and notations}
\end{figure}

Discretization of the governing equations is based on the semi-implicit method \cite{casulli1992} in which the barotropic pressure gradient and the vertical viscosity terms are discretized implicitly. However, the convective, the Coriolis, and the horizontal viscosity terms are finite differenced explicitly. A central difference scheme is used for spatial discretization. Notably, the nonlinear advection terms are solved based on the method described in \cite{shahabi2022} The general discretized form of the Eqs. (1) and (2) can be expressed as:

\begin{equation}
\begin{aligned}
    &\frac{u_{i-\frac{1}{2},j,k}^{n+1} - u_{i-\frac{1}{2},j,k}^n}{\Delta t} = F u_{i-\frac{1}{2},j,k}^n - g \frac{\zeta_{i,j}^{n+1} - \zeta_{i-1,j}^{n+1}}{\Delta x}+\frac{V_{i-\frac{1}{2},j,k-\frac{1}{2}} \frac{u_{i-\frac{1}{2},j,k-1}^{n+1} - u_{i-\frac{1}{2},j,k}^{n+1}}{\Delta z_{i-\frac{1}{2},j,k-\frac{1}{2}}} - V_{i-\frac{1}{2},j,k+\frac{1}{2}} \frac{u_{i-\frac{1}{2},j,k}^{n+1} - u_{i-\frac{1}{2},j,k+1}^{n+1}}{\Delta z_{i-\frac{1}{2},j,k+\frac{1}{2}}}}{\Delta z_{i-\frac{1}{2},j,k}} \\
\end{aligned}
\end{equation}
\begin{equation}
\begin{aligned}
    &\frac{v_{i,j-\frac{1}{2},k}^{n+1} - v_{i,j-\frac{1}{2},k}^n}{\Delta t} = F v_{i,j-\frac{1}{2},k}^n - g \frac{\zeta_{i,j}^{n+1} - \zeta_{i,j-1}^{n+1}}{\Delta y}+\frac{V_{i,j-\frac{1}{2},k-\frac{1}{2}} \frac{v_{i,j-\frac{1}{2},k-1}^{n+1} - v_{i,j-\frac{1}{2},k}^{n+1}}{\Delta z_{i,j-\frac{1}{2},k-\frac{1}{2}}} - V_{i,j-\frac{1}{2},k+\frac{1}{2}} \frac{v_{i,j-\frac{1}{2},k}^{n+1} - v_{i,j-\frac{1}{2},k+1}^{n+1}}{\Delta z_{i,j-\frac{1}{2},k+\frac{1}{2}}}}{\Delta z_{i,j-\frac{1}{2},k}},
\end{aligned}
\end{equation}

where \(\Delta t\) is the time step and \(F\) is an explicit operator including the explicit terms. There are several ways to discretize the explicit terms. We have used a second-order upwind explicit method for the convection terms and a central explicit method for the horizontal viscosity in the model.

By applying the dynamic boundary conditions, the discretized momentum equations for top layer cells take the form:

\begin{equation}
\begin{aligned}
    &\frac{u_{i-\frac{1}{2},j,1}^{n+1} - u_{i-\frac{1}{2},j,1}^n}{\Delta t} = F u_{i-\frac{1}{2},j,1}^n - g \frac{\zeta_{i,j}^{n+1} - \zeta_{i-1,j}^{n+1}}{\Delta x} + \frac{\tau_{wx} - V_{i-\frac{1}{2},j+\frac{1}{2}} \frac{u_{i-\frac{1}{2},j,1}^{n+1} - u_{i-\frac{1}{2},j,2}^{n+1}}{\Delta z_{i-\frac{1}{2},j+\frac{1}{2}}}}{\Delta z_{i-\frac{1}{2},j,1}},\\
\end{aligned}
\end{equation}
\begin{equation}
\begin{aligned}
    &\frac{v_{i,j-\frac{1}{2},1}^{n+1} - v_{i,j-\frac{1}{2},1}^n}{\Delta t} = F v_{i,j-\frac{1}{2},1}^n - g \frac{\zeta_{i,j}^{n+1} - \zeta_{i,j-1}^{n+1}}{\Delta y} + \frac{\tau_{wy} - V_{i,j-\frac{1}{2},1+\frac{1}{2}} \frac{v_{i,j-\frac{1}{2},1}^{n+1} - v_{i,j-\frac{1}{2},2}^{n+1}}{\Delta z_{i,j-\frac{1}{2},1+\frac{1}{2}}}}{\Delta z_{i,j-\frac{1}{2},1}},\\
\end{aligned}
\end{equation}

and for the bottom layer cells, the equations become:

\begin{equation}
\begin{aligned}
    &\frac{u_{i-\frac{1}{2},j,k}^{n+1} - u_{i-\frac{1}{2},j,k}^n}{\Delta t} = F u_{i-\frac{1}{2},j,k}^n - g \frac{\zeta_{i,j}^{n+1} - \zeta_{i-1,j}^{n+1}}{\Delta x} + \frac{V_{i-\frac{1}{2},j,k-\frac{1}{2}} \frac{u_{i-\frac{1}{2},j,k-1}^{n+1} - u_{i-\frac{1}{2},j,k}^{n+1}}{\Delta z_{i-\frac{1}{2},j,k-\frac{1}{2}}} - \tau_{bx}}{\Delta z_{i-\frac{1}{2},j,k}},\\
\end{aligned}
\end{equation}
\begin{equation}
\begin{aligned}
    &\frac{v_{i,j-\frac{1}{2},nk}^{n+1} - v_{i,j-\frac{1}{2},nk}^n}{\Delta t} = F v_{i,j-\frac{1}{2},nk}^n - g \frac{\zeta_{i,j}^{n+1} - \zeta_{i,j-1}^{n+1}}{\Delta y} + \frac{V_{i,j-\frac{1}{2},nk-\frac{1}{2}} \frac{v_{i,j-\frac{1}{2},nk-1}^{n+1} - v_{i,j-\frac{1}{2},nk}^{n+1}}{\Delta z_{i,j-\frac{1}{2},nk-\frac{1}{2}}} - \tau_{by}}{\Delta z_{i,j-\frac{1}{2},nk}}.\\
\end{aligned}
\end{equation}

where the velocities in the time step \(n\) are used to calculate \(\tau_{bx}\) and \(\tau_{by}\).

A detailed description of the algorithm is presented below. For brevity, the \(i\) and \(j\) indices are not shown if they are unnecessary. The Eqs. (12)-(15) can be written in the following form:

\begin{equation}
\begin{aligned}
    &\mathbf{A}_{i-\frac{1}{2},j}^n \mathbf{U}_{i-\frac{1}{2},j}^{n+1} = \mathbf{D}_{i-\frac{1}{2},j}^{n+1}, \\
\end{aligned}
\end{equation}
\begin{equation}
\begin{aligned}
    &\mathbf{A}_{i,j-\frac{1}{2}}^n \mathbf{V}_{i,j-\frac{1}{2}}^{n+1} = \mathbf{D}_{i,j-\frac{1}{2}}^{n+1},
\end{aligned}
\end{equation}

where

\begin{equation*}
\begin{aligned}
    &\mathbf{A}^n =
\begin{bmatrix}
b_1 & c_1 & & & \\
a_2 & b_2 & c_2 & & \\
& \ddots & \ddots & \ddots & \\
& & a_{nk-1} & b_{nk-1} & c_{nk-1} \\
& & & a_{nk} & b_{nk}
\end{bmatrix},
\quad
\mathbf{U}^{n+1} =
\begin{bmatrix}
u_1^{n+1} \\
u_2^{n+1} \\
\vdots \\
u_{nk-1}^{n+1} \\
u_{nk}^{n+1}
\end{bmatrix},
\quad
\mathbf{V}^{n+1} =
\begin{bmatrix}
v_1^{n+1} \\
v_2^{n+1} \\
\vdots \\
v_{nk-1}^{n+1} \\
v_{nk}^{n+1}
\end{bmatrix},
\quad
\mathbf{D}^{n+1} =
\begin{bmatrix}
d_1 \\
d_2 \\
\vdots \\
d_{nk-1} \\
d_{nk}
\end{bmatrix}.
\end{aligned}
\end{equation*}

wherein the coefficients are:
\begin{equation}
\begin{aligned}
    &a_k = -\frac{v_{z_{k-\frac{1}{2}}} \Delta t}{\Delta z_{k-\frac{1}{2}} \Delta z_k},\\
    &c_k = \frac{v_{z_{k+\frac{1}{2}}} \Delta t}{\Delta z_{k+\frac{1}{2}} \Delta z_k},\\
    &b_k = 1 - (a_k + c_k),\\
    &d_{i-\frac{1}{2},j,k} = d'_{i-\frac{1}{2},j,k} - g \Delta t \frac{\zeta_{i,j}^{n+1} - \zeta_{i-1,j}^{n+1}}{\Delta x},\\
    &d_{i,j-\frac{1}{2},k} = d'_{i,j-\frac{1}{2},k} - g \Delta t \frac{\zeta_{i,j}^{n+1} - \zeta_{i,j-1}^{n+1}}{\Delta y},\\
    &d'_{i-\frac{1}{2},j,k} = u_{i-\frac{1}{2},j,k}^n + \Delta t F u_{i-\frac{1}{2},j,k}^n,\\
    &d'_{i,j-\frac{1}{2},k} = v_{i,j-\frac{1}{2},k}^n + \Delta t F v_{i,j-\frac{1}{2},k}^n.
    \end{aligned}
\end{equation}

Now applying the dynamic boundary conditions yields to:

\begin{equation}
\begin{aligned}
    &a_1 = c_{nk} = 0,\\
&d'_{i-\frac{1}{2},j,1} = u_{i-\frac{1}{2},j,1}^n + F \Delta t u_{i-\frac{1}{2},j,1}^n + \frac{\Delta t \tau_{wx}}{\Delta z_{i-\frac{1}{2},j,1}},\\
&d'_{i-\frac{1}{2},j,nk} = u_{i-\frac{1}{2},j,nk}^n + F \Delta t u_{i-\frac{1}{2},j,nk}^n - \frac{\Delta t \tau_{bx}}{\Delta z_{i-\frac{1}{2},j,nk}},\\
&d'_{i,j-\frac{1}{2},1} = v_{i,j-\frac{1}{2},1}^n + F \Delta t v_{i,j-\frac{1}{2},1}^n + \frac{\Delta t \tau_{wy}}{\Delta z_{i,j-\frac{1}{2},1}},\\
&d'_{i,j-\frac{1}{2},nk} = v_{i,j-\frac{1}{2},nk}^n + F \Delta t v_{i,j-\frac{1}{2},nk}^n + \frac{\Delta t \tau_{by}}{\Delta z_{i,j-\frac{1}{2},nk}}.
\end{aligned}
\end{equation}

The proposed algorithm for solving Eqs. (16-17) consists of three steps. First, a forward sweep eliminates the lower diagonal of the coefficients matrix. Second, the right-hand side entries are partitioned into the known and unknown parts and each part is formulated recursively. Third, a backward sweep solves the system for \(U^{n+1}\) and \(V^{n+1}\).

\subsection{Step 1: Forward Sweep}
In the first step, the lower diagonal is eliminated by a forward sweep. This step is similar to the well-known double-sweep algorithm \cite{anderson1995}. By applying the forward sweep, the Eqs. (16) and (17) become as follows:

\begin{equation}
\begin{aligned}
    &\mathbf{E}_{i-\frac{1}{2},j}^n \mathbf{U}_{i-\frac{1}{2},j}^{n+1} = \mathbf{F}_{i-\frac{1}{2},j}^{n+1},\\
\end{aligned}
\end{equation}
\begin{equation}
\begin{aligned}
    &\mathbf{E}_{i,j-\frac{1}{2}}^n \mathbf{V}_{i,j-\frac{1}{2}}^{n+1} = \mathbf{F}_{i,j-\frac{1}{2}}^{n+1},\\
\end{aligned}
\end{equation}

where

\begin{equation*}
\begin{aligned}
    &\mathbf{E}^n =
    \begin{bmatrix}
    1 & -e_1^n & & & \\
    & \ddots & \ddots & & \\
    & & 1 & -e_k^n & \\
    & & & \ddots & \\
    & & & & 1
    \end{bmatrix},
    \quad
    \mathbf{F}^{n+1} =
    \begin{bmatrix}
    f_1^{n+1} \\
    \vdots \\
    f_k^{n+1} \\
    \vdots \\
    f_{nk}^{n+1}
    \end{bmatrix}.
    \end{aligned}
\end{equation*}

wherein the entries are:
\begin{equation}
\begin{aligned}
    &e_k^n = \frac{-c_k}{a_k e_{k-1}^n + b_k}, \quad e_0^n = 0,\\
\end{aligned}
\end{equation}
\begin{equation}
\begin{aligned}
    &f_k^{n+1} = \frac{d_k^{n+1} - a_k f_{k-1}^{n+1}}{a_k e_{k-1}^n + b_k}, \quad f_0^{n+1} = 0.\\
\end{aligned}
\end{equation}

\subsection{Step 2: Decomposing \(F^{n+1}\)}
In the second step, the right-hand side matrix entries, i.e. Eq. (23), are partitioned into two parts in a way that each part can be defined recursively.

\textbf{Lemma 1.} The Eq. (23) can be defined as follows for the \(x\) and \(y\) directions respectively:

\begin{equation}
\begin{aligned}
    &f_{i-\frac{1}{2},j,k}^{n+1} = \psi_{i-\frac{1}{2},j,k}^{(1)} + \psi_{i-\frac{1}{2},j,k}^{(2)} \left( \zeta_{i,j}^{n+1} - \zeta_{i-1,j}^{n+1} \right), \\
\end{aligned}
\end{equation}
\begin{equation}
\begin{aligned}
    &f_{i,j-\frac{1}{2},k}^{n+1} = \psi_{i,j-\frac{1}{2},k}^{(1)} + \psi_{i,j-\frac{1}{2},k}^{(2)} \left( \zeta_{i,j}^{n+1} - \zeta_{i,j-1}^{n+1} \right),
\end{aligned}
\end{equation}

where the recursive coefficient \(\psi^{(1)}\) for both \(x\) and \(y\) directions is:

\begin{equation}
\begin{aligned}
    \psi_k^{(1)} = \frac{d_k' - a_k \psi_{k-1}^{(1)}}{a_k e_{k-1}^n + b_k}, \quad \psi_0^{(1)} = 0,
\end{aligned}
\end{equation}

and \(\psi^{(2)}\) for \(x\) and \(y\) directions are:

\begin{equation}
\begin{aligned}
    &\psi_{i-\frac{1}{2},j,k}^{(2)} = \frac{-g \frac{\Delta t}{\Delta x} - a_k \psi_{i-\frac{1}{2},j,k-1}^{(2)}}{a_k e_{k-1}^n + b_k}, \quad \psi_{i-\frac{1}{2},j,0}^{(2)} = 0,\\
\end{aligned}
\tag{26b}
\end{equation}
\begin{equation}
\begin{aligned}
    &\psi_{i,j-\frac{1}{2},k}^{(2)} = \frac{-g \frac{\Delta t}{\Delta y} - a_k \psi_{i,j-\frac{1}{2},k-1}^{(2)}}{a_k e_{k-1}^n + b_k}, \quad \psi_{i,j-\frac{1}{2},0}^{(2)} = 0.
\end{aligned}
\tag{26c}
\end{equation}

\textbf{Proof.} Using mathematical induction, the Eq. (24) is proved and likewise the Eq. (25) can be proved. In the basis step, it is proved that the Eq. (24) holds for \(k=1\). From the Eq. (23), \(f_1\) is:

\begin{equation*}
\begin{aligned}
    &f_1^{n+1} = \frac{d_1 - a_1 f_0^{n+1}}{a_1 e_0 + b_1} = \frac{\left( d_1' - g \frac{\Delta t}{\Delta x} (\zeta_{i,j}^{n+1} - \zeta_{i-1,j}^{n+1}) \right) - 0}{0 + b_1} = \frac{d_1'}{b_1} - \frac{g \Delta t}{b_1 \Delta x} \left( \zeta_{i,j}^{n+1} - \zeta_{i-1,j}^{n+1} \right),\\
    &\Rightarrow f_1^{n+1} = \psi_1^{(1)} + \psi_1^{(2)} (\zeta_{i,j}^{n+1} - \zeta_{i-1,j}^{n+1}).
\end{aligned}
\end{equation*}

In the inductive step, it is proved if the Eq. (24) holds for \(k\), then it holds for \(k+1\). We have:

\begin{equation*}
\begin{aligned}
    &f_{k+1}^{n+1} = \frac{d_{k+1} - a_{k+1} f_{k}^{n+1}}{a_{k+1} e_{k}^n + b_{k+1}} \quad (*) \\
    &f_{k}^{n+1} = \psi_{k}^{(1)} + \psi_{k}^{(2)} (\zeta_{i,j}^{n+1} - \zeta_{i-1,j}^{n+1}) \quad (**).
\end{aligned}
\end{equation*}

Substituting (**) into (*) gives:

\begin{equation*}
\begin{aligned}
    &f_{k+1}^{n+1} = \frac{\left( d_{k+1}' - g \frac{\Delta t}{\Delta x} (\zeta_{i,j}^{n+1} - \zeta_{i-1,j}^{n+1}) \right) - a_{k+1} \left( \psi_{k}^{(1)} + \psi_{k}^{(2)} (\zeta_{i,j}^{n+1} - \zeta_{i-1,j}^{n+1}) \right)}{a_{k+1} e_{k}^n + b_{k+1}}\\
    &= \left( \frac{d_{k+1}' - a_{k+1} \psi_{k}^{(1)}}{a_{k+1} e_{k}^n + b_{k+1}} \right) + \left( \frac{-g \frac{\Delta t}{\Delta x} - a_{k+1} \psi_{k}^{(2)}}{a_{k+1} e_{k}^n + b_{k+1}} \right) (\zeta_{i,j}^{n+1} - \zeta_{i-1,j}^{n+1})\\
    &\Rightarrow f_{k+1}^{n+1} = \psi_{k+1}^{(1)} + \psi_{k+1}^{(2)} (\zeta_{i,j}^{n+1} - \zeta_{i-1,j}^{n+1}).
\end{aligned}
\end{equation*}

\subsection{Step 3: Backward Sweep}
In the third step, a backward sweep produces the horizontal velocities in terms of the free-surface elevations.

\textbf{Lemma 2.} The horizontal velocity components can be expressed as follows:

\begin{equation}
\begin{aligned}
    &u_{i-\frac{1}{2},j,k}^{n+1} = \Omega_{i-\frac{1}{2},j,k}^{(1)} + \Omega_{i-\frac{1}{2},j,k}^{(2)} (\zeta_{i,j}^{n+1} - \zeta_{i-1,j}^{n+1}), \\
\end{aligned}
\end{equation}
\begin{equation}
\begin{aligned}
    &v_{i,j-\frac{1}{2},k}^{n+1} = \Omega_{i,j-\frac{1}{2},k}^{(1)} + \Omega_{i,j-\frac{1}{2},k}^{(2)} (\zeta_{i,j}^{n+1} - \zeta_{i,j-1}^{n+1}),\\
\end{aligned}
\end{equation}

where the coefficients \(\Omega^{(1)}\) and \(\Omega^{(2)}\) are defined as:

\begin{equation}
\begin{aligned}
    &\Omega_k^{(1)} = \psi_k^{(1)} + e_k^n \Omega_{k+1}^{(1)}, \quad \Omega_{nk+1}^{(1)} = 0, \\
\end{aligned}
\end{equation}
\begin{equation}
\begin{aligned}
    &\Omega_k^{(2)} = \psi_k^{(2)} + e_k^n \Omega_{k+1}^{(2)}, \quad \Omega_{nk+1}^{(2)} = 0.\\
\end{aligned}
\tag{29b}
\end{equation}

\textbf{Proof.} A backward induction is used to prove Eqs. (27) and (28). In the basis step, it is proved that the Eq. (27) holds for \(k=n_k\). The Eq. (20) for \(k=n_k\) is:

\begin{equation*}
\begin{aligned}
    &u_{nk}^{n+1} = f_{nk}^{n+1}\\
    &= \psi_{nk}^{(1)} + \psi_{nk}^{(2)} (\zeta_{i,j}^{n+1} - \zeta_{i-1,j}^{n+1})\\    
    &\Rightarrow u_{nk}^{n+1} = \Omega_{nk}^{(1)} + \Omega_{nk}^{(2)} (\zeta_{i,j}^{n+1} - \zeta_{i-1,j}^{n+1}).\\
\end{aligned}
\end{equation*}

In the inductive step, it is shown if the Eq. (27) holds for \(k+1\), it holds for \(k\). the \(k\)-th row of the Eq. (20) is:

\begin{equation*}
\begin{aligned}
    &u_k^{n+1} - e_k u_k^{n+1} = f_k^{n+1} \quad (*).
\end{aligned}
\end{equation*}
and Eq. (27) holds for \(k+1\) is:
\begin{equation*}
\begin{aligned}
    &u_{k+1}^{n+1} = \Omega_{k+1}^{(1)} + \Omega_{k+1}^{(2)} (\zeta_{i,j}^{n+1} - \zeta_{i-1,j}^{n+1}) \quad (**).
\end{aligned}
\end{equation*}

Substituting (**) into (*) yields:

\begin{equation*}
\begin{aligned}
    &u_k^{n+1} - e_k \left( \Omega_{k+1}^{(1)} + \Omega_{k+1}^{(2)} (\zeta_{i,j}^{n+1} - \zeta_{i-1,j}^{n+1}) \right) = \psi_k^{(1)} + \psi_k^{(2)} (\zeta_{i,j}^{n+1} - \zeta_{i-1,j}^{n+1})\\
    &\Rightarrow u_k^{n+1} = \left( \psi_k^{(1)} + e_k^n \Omega_{k+1}^{(1)} \right) + \left( \psi_k^{(2)} + e_k^n \Omega_{k+1}^{(2)} \right) (\zeta_{i,j}^{n+1} - \zeta_{i-1,j}^{n+1})\\
    &\Rightarrow u_k^{n+1} = \Omega_k^{(1)} + \Omega_k^{(2)} (\zeta_{i,j}^{n+1} - \zeta_{i-1,j}^{n+1}).
\end{aligned}
\end{equation*}

The Eq. (28) can be also proved similarly. Finally, the horizontal velocity components are determined in terms of free-surface elevations as follows (matrix form):

\begin{equation}
\begin{aligned}
    &\mathbf{U}_{i-\frac{1}{2},j}^{n+1} = \mathbf{\Omega}_{i-\frac{1}{2},j}^{(1)} + \mathbf{\Omega}_{i-\frac{1}{2},j}^{(2)} (\zeta_{i,j}^{n+1} - \zeta_{i-1,j}^{n+1}).\\
\end{aligned}
\end{equation}
\begin{equation}
\begin{aligned}
    &\mathbf{V}_{i,j-\frac{1}{2}}^{n+1} = \mathbf{\Omega}_{i,j-\frac{1}{2}}^{(1)} + \mathbf{\Omega}_{i,j-\frac{1}{2}}^{(2)} (\zeta_{i,j}^{n+1} - \zeta_{i,j-1}^{n+1}).\\
\end{aligned}
\end{equation}

In order to determine \(\xi^{n+1}\), the Eqs. (30) and (31) are substituted into the depth-integrated continuity equation. The discretized form of the Eq. (6) in a matrix notation takes the form:

\begin{equation}
\begin{aligned}
    &\frac{\zeta_{i,j}^{n+1} - \zeta_{i,j}^n}{\Delta t} + \frac{(\mathbf{\Delta z})_{i+\frac{1}{2},j}^{\top} \mathbf{U}_{i+\frac{1}{2},j}^{n+1} - (\mathbf{\Delta z})_{i-\frac{1}{2},j}^{\top} \mathbf{U}_{i-\frac{1}{2},j}^{n+1}}{\Delta x} + \frac{(\mathbf{\Delta z})_{i,j+\frac{1}{2}}^{\top} \mathbf{U}_{i,j+\frac{1}{2}}^{n+1} - (\mathbf{\Delta z})_{i,j-\frac{1}{2}}^{\top} \mathbf{U}_{i,j-\frac{1}{2}}^{n+1}}{\Delta y} = 0.\\
\end{aligned}
\end{equation}

where \(A\) is a five-diagonal matrix and \(B\) is the right-hand side vector.

Substitution of the Eqs. (30) and (31) into the Eq. (32) yields:

\begin{equation}
\begin{aligned}
    &\zeta_{i,j}^{n+1} + \frac{\Delta t}{\Delta x} \left\{ \left[ (\mathbf{\Delta z})_{i+\frac{1}{2},j}^{\top} \mathbf{\Omega}_{i+\frac{1}{2},j}^{(2)} \right] (\zeta_{i+1,j}^{n+1} - \zeta_{i,j}^{n+1}) - \left[ (\mathbf{\Delta z})_{i-\frac{1}{2},j}^{\top} \mathbf{\Omega}_{i-\frac{1}{2},j}^{(2)} \right] (\zeta_{i,j}^{n+1} - \zeta_{i-1,j}^{n+1}) \right\}\\
    &+ \frac{\Delta t}{\Delta y} \left\{ \left[ (\mathbf{\Delta z})_{i,j+\frac{1}{2}}^{\top} \mathbf{\Omega}_{i,j+\frac{1}{2}}^{(2)} \right] (\zeta_{i,j+1}^{n+1} - \zeta_{i,j}^{n+1}) - \left[ (\mathbf{\Delta z})_{i,j-\frac{1}{2}}^{\top} \mathbf{\Omega}_{i,j-\frac{1}{2}}^{(2)} \right] (\zeta_{i,j}^{n+1} - \zeta_{i,j-1}^{n+1}) \right\}\\
    &= \zeta_{i,j}^n - \frac{\Delta t}{\Delta x} \left\{ \left[ (\mathbf{\Delta z})_{i+\frac{1}{2},j}^{\top} \mathbf{\Omega}_{i+\frac{1}{2},j}^{(1)} \right] - \left[ (\mathbf{\Delta z})_{i-\frac{1}{2},j}^{\top} \mathbf{\Omega}_{i-\frac{1}{2},j}^{(1)} \right] \right\}\\
    &- \frac{\Delta t}{\Delta y} \left\{ \left[ (\mathbf{\Delta z})_{i,j+\frac{1}{2}}^{\top} \mathbf{\Omega}_{i,j+\frac{1}{2}}^{(1)} \right] - \left[ (\mathbf{\Delta z})_{i,j-\frac{1}{2}}^{\top} \mathbf{\Omega}_{i,j-\frac{1}{2}}^{(1)} \right] \right\},\\
\end{aligned}
\end{equation}

which is a symmetric positive definite and diagonally dominant five-diagonal system. This system can be solved efficiently by a conjugate gradient method \cite{casulli1992}.

In summary, the proposed algorithm to calculate \(\xi^{n+1}\) takes the following simple steps:
\begin{enumerate}
    \item Calculate \(a_k, b_k, c_k\) and \(d'_k\) from the Eqs. (21) and (22).
    \item Calculate \(e_k\) from the Eq. (25).
    \item Calculate \(\psi^{(1)}_k\) and \(\psi^{(2)}_k\) from the Eqs. (29).
    \item Calculate \(\Omega^{(1)}_k\) and \(\Omega^{(2)}_k\) from the Eqs. (32).
    \item Solve the five-diagonal linear system (36).
\end{enumerate}

\subsection{Calculating Velocity Field}
Once the new free-surface elevations have been determined, the new horizontal velocity components are directly determined by substituting the new free-surface elevations in the Eqs. (27) and (28). Eventually, the Eq. (3) is directly applicable to calculate the vertical velocity components. This equation is discretized as follows:

\begin{equation}
\begin{aligned}
    &\frac{u_{i+\frac{1}{2},j,k}^{n+1} - u_{i-\frac{1}{2},j,k}^{n+1}}{\Delta x} + \frac{v_{i,j+\frac{1}{2},k}^{n+1} - v_{i,j-\frac{1}{2},k}^{n+1}}{\Delta y} + \frac{w_{i,j,k+\frac{1}{2}}^{n+1} - w_{i,j,k-\frac{1}{2}}^{n+1}}{\Delta z_{i,k}} = 0.\\
\end{aligned}
\end{equation}

Using the boundary condition (4) that implies \(w^{n+1}_{n_k}=0\) and starting from the bottom layer cells, the vertical velocity components are determined.

\section{Method Validation and Efficiency}
Two classical test cases have been undertaken to validate the proposed algorithm and to verify the mass and momentum conservation of the developed model.

\subsection{Standing Wave in a Rectangular Basin}
In the first test case, a standing wave in a rectangular basin was simulated. The analytical solution for a small amplitude standing wave is given by the following relations \cite{mei1989}:

\begin{equation}
    \xi = A \cos(k_x x) \cos(k_y y) \cos(\omega t),\\
\end{equation}
\begin{equation}
    (u, v) = A \frac{g}{\sigma} \frac{\cosh k \left( h + (-z) \right)}{\cosh kh} \sin(\sigma t) \left( \left( k_x \sin(k_x x) \cos(k_y y), k_y \sin(k_x x) \cos(k_y y) \right) \right),\\
\end{equation}
\begin{equation}
    w = A \frac{gk}{\sigma} \frac{\sinh k \left( h + (-z) \right)}{\cosh kh} \cos(k_x x) \cos(k_y y) \sin(\sigma t),\\
\end{equation}

where $A$ is the wave amplitude, $\left(k_x, k_y\right) = \left(\frac{\pi}{L_x}, \frac{\pi}{L_y}\right)$ which in this test $L_x$ and $L_y$ are taken as the basin dimensions, $k = \sqrt{k_x^2 + k_y^2}$, $h$ is the mean water depth, $z$ is the vertical position, and $\sigma = k \sqrt{gh}$. Figure 2 shows the initial condition for the free-surface elevation. The analytical solution is derived in an ideal condition that the fluid is inviscid and the bottom friction is zero. These conditions have been also applied in the model. The wave amplitude and the basin dimensions have been taken as 0.1m and $500 \times 500 \times 10$m. The computational space has been discretized with the cells of dimensions $\Delta x = 10$m, $\Delta y = 10$m, $\Delta z = 1$m, and the time step was 0.05s.

\begin{figure}[ht]
\centering
\includegraphics[width=0.8\textwidth]{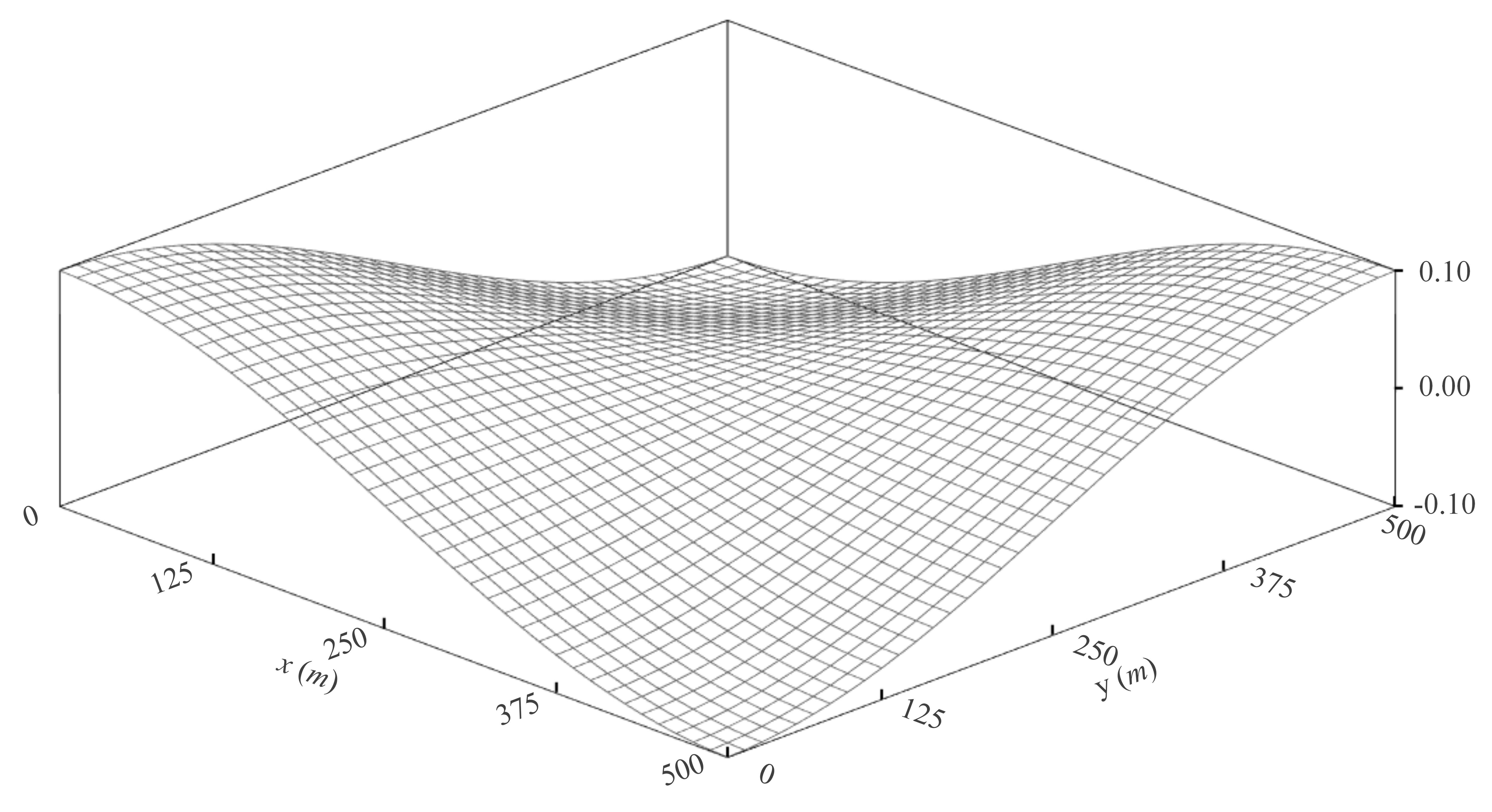}
\caption{The initial condition of the free-surface in the 3D standing wave test case}
\end{figure}

Figure 3(a) shows the numerical prediction of the free-surface elevations compared with the analytical solution at \(t=6T\). As shown in this figure, wave damping after six oscillations is small and acceptable. In figure 3(b), numerical and analytical solutions of the free-surface elevation are presented at \(x=0, y=0\). Figure 3(c) shows the predicted and analytical values of the horizontal velocity \(u\) at \(x=250\)m, \(y=0\). It can be seen that the simulated free-surface and horizontal velocity agree well with the analytical solution.

\begin{figure}[ht]
\centering
\includegraphics[width=0.8\textwidth]{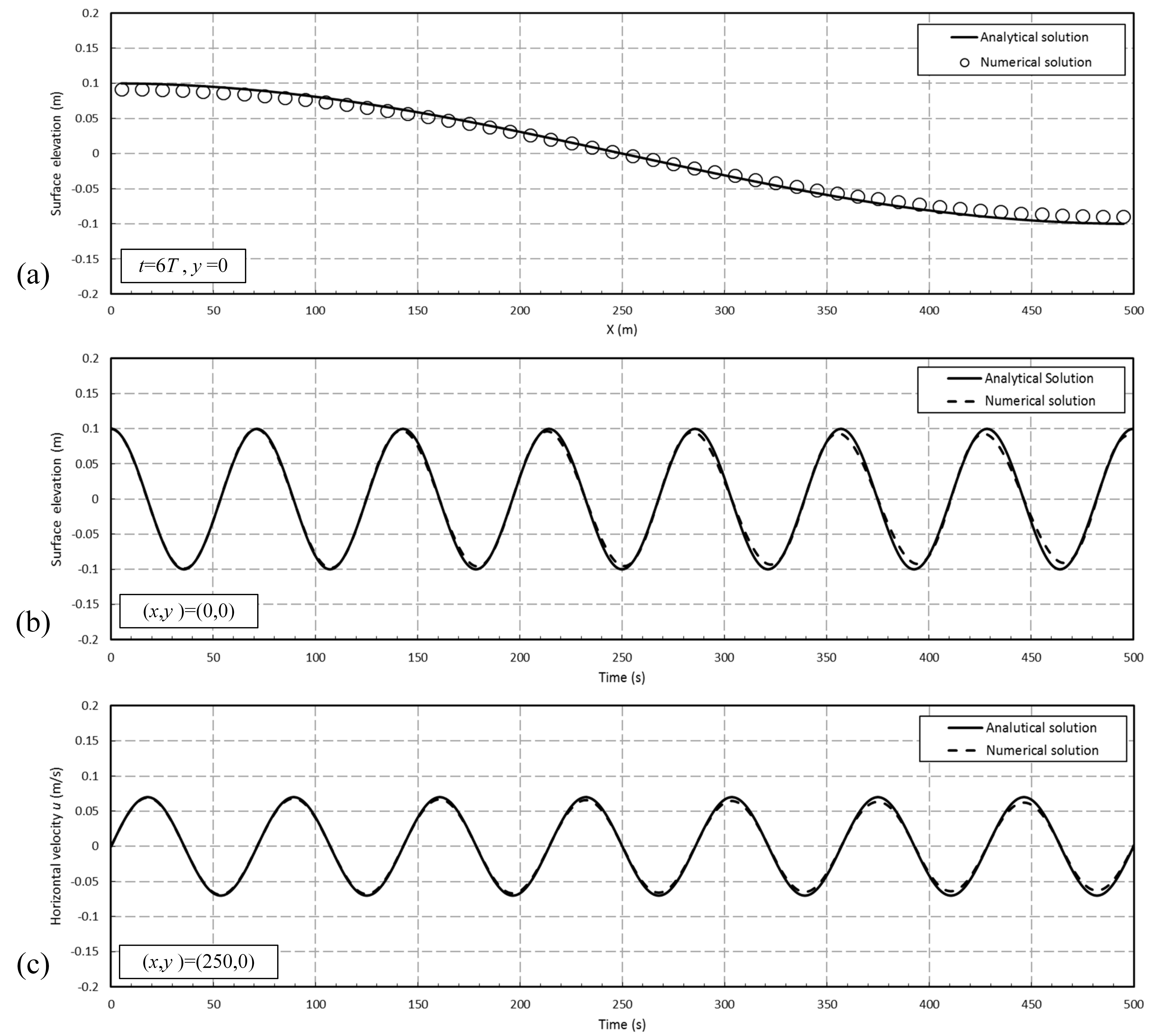}
\caption{Comparisons of simulated free-surface and horizontal velocity \(u\) with the analytical solution for the 3D linear standing wave test case: (a) free-surface elevations at \(y=0\) and \(t=6T\); (b) free-surface elevation at \((x, y)=(0, 0)\); (c) horizontal velocity \(u\) at \((x, y)=(250, 0)\).}
\end{figure}

\subsection{Wind-Driven Circulatory Flow}
The second test case was a wind-driven circulatory flow in a closed rectangular basin. In this test, the advection and the horizontal diffusion were neglected and a linearized bottom friction was used. The analytical solution for a constant vertical eddy viscosity and a linearized bottom friction is \cite{huang1995}:

\begin{equation}
\begin{aligned}
    &u = \frac{1}{6 \nu_z} g \frac{\partial \xi}{\partial \tilde{x}} \left( 3z^2 - H^2 \right) + \frac{\tau_w}{2 \rho \nu_z} \left( H + 2z \right),\\
\end{aligned}
\end{equation}
\begin{equation}
\begin{aligned}
    &\frac{\partial \xi}{\partial x} = \frac{3}{2} \frac{\tau_w}{\rho g H} \frac{(2 \nu_z + k H)}{(3 \nu_z + k H)},\\
\end{aligned}
\end{equation}

where \(H\) is the water depth and \(k\) is the linearized bottom friction coefficient. In this test, \(\tau_w\) was set to 0.1 N/m\(^2\) in both numerical and analytical solutions. The other parameters were: \(H=40\)m, \(\nu_z=0.03\)m\(^2\)/s, \(k=0.005\)m/s, and \(\rho=1000\)Kg/m\(^3\). The basin dimensions were 2500×2500×40m, which were discretized into the cells with \(\Delta x=50\) m, \(\Delta x=50\) m, \(\Delta z=2\) m, and the time step was 2s. A comparison of the horizontal velocity profile between the numerical model and the analytical solution is shown in Figure 4.

\begin{figure}[ht]
\centering
\includegraphics[width=0.6\textwidth]{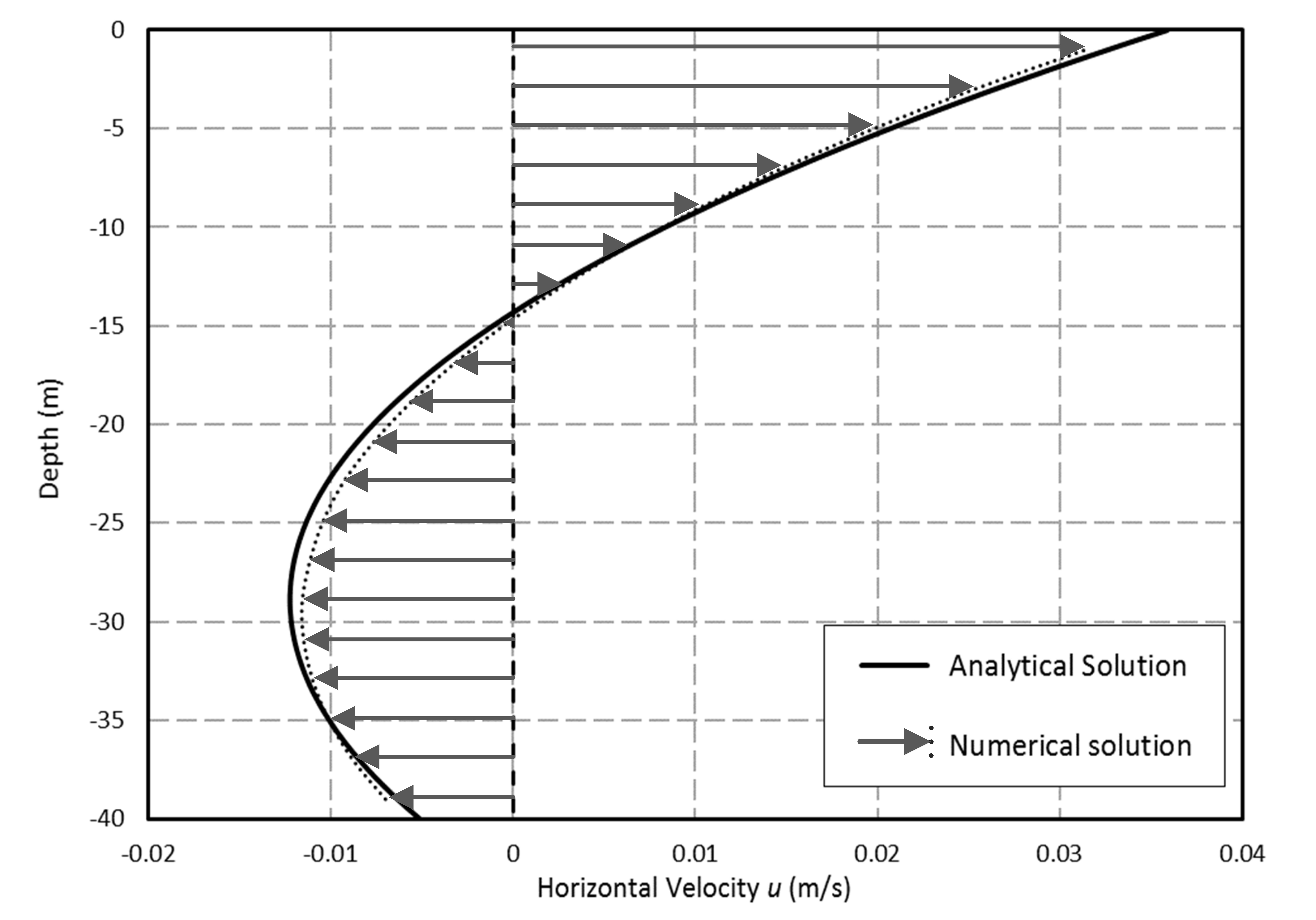}
\caption{Comparison of the simulated horizontal velocity profile with the analytical solution at \(x=1250\)m, \(y=1250\)m for the wind-driven circulatory flow case}
\end{figure}

\subsection{Efficiency of the Method}
Assuming that the number of horizontal layers is \(N\), any methods employed to couple the momentum and the continuity equations will finish at \(O(N)\), \(O(N^2)\), etc. Generally, the smaller the order of complexity of the algorithm, the faster it will run. The standing wave test case in section 4.1 was repeated for a different number of horizontal layers to investigate the order of complexity of the method and all CPU times were recorded.

Figure 5(a) shows the simulation time for a different number of horizontal layers. As shown in this figure, the simulation time increases linearly and the method is of the least order of complexity \(O(N)\).

\begin{figure}[ht]
\centering
\includegraphics[width=0.6\textwidth]{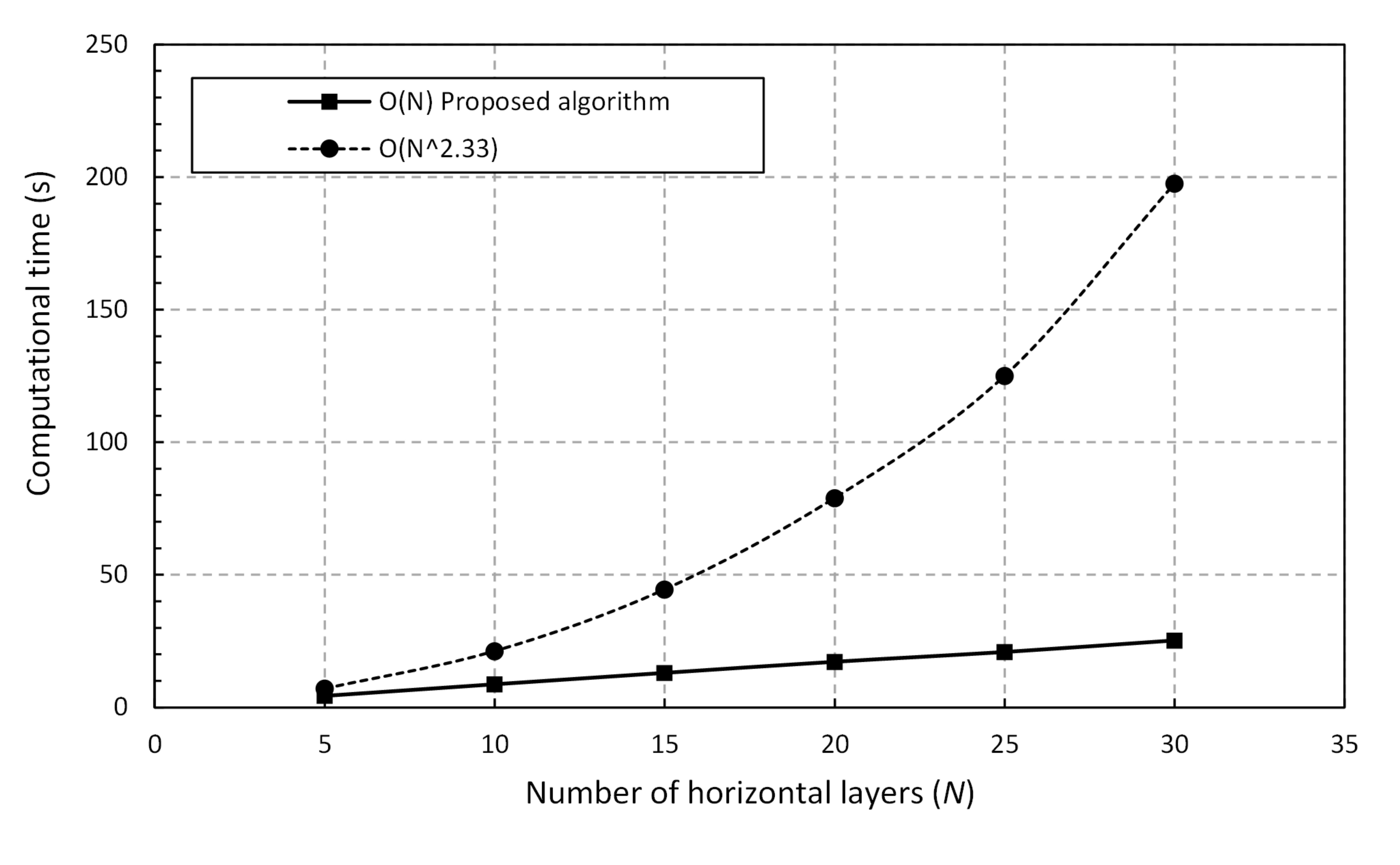}
\caption{Computation time versus the number of horizontal layers}
\end{figure}

For the purpose of comparison, a fast tridiagonal matrix inversion algorithm \cite{usmani1994} with the order of complexity \(O(N^{2.33})\) was also used to couple the equations in the test case. The result shows how important is the order of complexity of the employed coupling method in the computational efficiency of the model.

\section{Summary}
A semi-implicit model is developed to simulate three-dimensional HPEs. The pressure gradient terms and the vertical viscosity terms are discretized implicitly in the horizontal momentum equations, while the other terms are treated explicitly with respect to time. A method is presented to couple the horizontal momentum equation and the depth-integrated continuity equation which constructs a linear system to calculate the free-surface elevations. In the proposed algorithm, some simple recursive formulas have been formulated to express the horizontal velocity components in terms of the neighbor free-surface elevations which allows them to be directly substituted in the depth-integrated continuity equation.

The algorithm was validated by two idealized classical test cases: a uni-nodal standing wave in a rectangular basin and a wind-driven circulatory flow in a closed rectangular basin. It has been shown that the computation time increases linearly when the number of horizontal layers increases. It means that the order of complexity of the method is \(O(N)\). Although the most expensive part of computations in each time step is to solve the five-diagonal linear system of the free-surface, it has been shown that the order of complexity the employed method to couple the equations considerably affects the computational efficiency of the model.

\bibliographystyle{unsrt}  


\begin{thebibliography}{9}

\bibitem{anderson1995}
Anderson J. (1995). Computational Fluid Dynamics (Vol. 206). New York: McGraw-Hill. \doi{https://doi.org/10.1201/9781351124027}.

\bibitem{blumberg1987}
Blumberg A. F., \& Mellor G. L. (1987). A description of a three‐dimensional coastal ocean circulation model. Three‐dimensional coastal ocean models, 4, 1-16. \doi{https://doi.org/10.1029/CO004p0001}.

\bibitem{blumberg1994}
Blumberg A. F. (1994). A primer for ECOM-si. Technical report of HydroQual 66.

\bibitem{casulli1992}
Casulli V., \& Cheng R. T. (1992). Semi-implicit finite difference methods for three-dimensional shallow water flow. International Journal for Numerical Methods in Fluids, 15(6), 629-648. \doi{https://doi.org/10.1002/fld.1650150602}.

\bibitem{casulli1997}
Casulli V. (1997). Numerical simulation of three-dimensional free surface flow in isopycnal coordinates. International Journal for Numerical Methods in Fluids, 25(6), 645-658. \doi{https://doi.org/10.1002/(SICI)1097-0363(19970930)25:6\%3C645::AID-FLD579\%3E3.0.CO;2-L}.

\bibitem{casulli1998}
Casulli V., \& Stelling G. S. (1998). Numerical Simulation of 3D Quasi-Hydrostatic Free-Surface Flows. Journal of Hydraulic Engineering, 124(7), 678-686. \doi{https://doi.org/10.1061/(ASCE)0733-9429(1998)124:7(678)}.

\bibitem{casulli1999}
Casulli V. (1999). A Semi-Implicit Finite Difference Method for Non-Hydrostatic Free-Surface Flows. International Journal for Numerical Methods in Fluids, 30(May 1998), 425-440. \doi{https://doi.org/10.1002/(SICI)1097-0363(19990630)30:4\%3C425::AID-FLD847\%3E3.0.CO;2-D}.

\bibitem{casulli2000}
Casulli V., \& Walters R. A. (2000). An unstructured grid three-dimensional model based on the shallow water equations. International Journal for Numerical Methods in Fluids, 32(3), 331-348. \doi{https://doi.org/10.1002/(SICI)1097-0363(20000215)32:3\%3C331::AID-FLD941\%3E3.0.CO;2-C}.

\bibitem{chen2003}
Chen C., Liu H., \& Beardsley R. C. (2003). An unstructured finite volume three-dimensional primitive equation ocean model: Application to coastal ocean and estuaries. Journal Atmospheric and Oceanic Technology, 20, 159-186. \doi{https://doi.org/10.1175/1520-0426(2003)020\%3C0159:AUGFVT\%3E2.0.CO;2}.

\bibitem{chen2003a}
Chen X. (2003). A free-surface correction method for simulating shallow water flows. Journal of Computational Physics, 189(2), 557-578. \doi{https://doi.org/10.1016/S0021-9991(03)00234-1}.

\bibitem{fringer2006}
Fringer O. B., Gerritsen M., \& Street R. L. (2006). An unstructured-grid finite-volume nonhydrostatic parallel coastal ocean simulator. Ocean Modelling, 14(3-4), 139-173. \doi{https://doi.org/10.1016/j.ocemod.2006.03.006}.

\bibitem{huang1995}
Huang, H. \& Spaulding, M. (1995). A three-dimensional numerical model of estuarine circulation and water quality induced by surface discharges. Estuarine, Coastal and Shelf Science, 40, 357-380.

\bibitem{madala1977}
Madala R. V., \& Piacseki S. A. (1977). A semi-implicit numerical model for baroclinic oceans. Journal of Computational Physics, 23(2), 167-178. \doi{https://doi.org/10.1016/0021-9991(77)90119-X}.

\bibitem{marshall1997}
Marshall J., Hill C., Perelman L., \& Adcroft A. (1997). Hydrostatic, quasi-hydrostatic, and nonhydrostatic ocean modeling. Journal of Geophysical Research: Oceans, 102(C3), 5733-5752. \doi{https://doi.org/10.1029/96JC02776}.

\bibitem{mei1989}
Mei C. C. (1989). The applied dynamics of ocean surface waves. Singapore: World Scientific.

\bibitem{shahabi2022}
Shahabi, A. and Ghiassi, R. (2022). A robust second-order godunov-type method for Burgers’ equation. International Journal of Applied and Computational Mathematics, 8(2), 82. Springer. \doi{https://doi.org/10.1007/s40819-021-01171-7}.

\bibitem{shchepetkin2005}
Shchepetkin A. F., \& McWilliams J. C. (2005). The regional oceanic modeling system (ROMS): a split-explicit, free-surface, topography-following-coordinate oceanic model. Ocean Modelling, 9(4), 347-404. \doi{https://doi.org/10.1016/j.ocemod.2004.08.002}.

\bibitem{smagorinsky1963}
Smagorinsky J. (1963). General Circulation Experiments With the Primitive Equations. Monthly Weather Review, 91(3), 99-164.

\bibitem{taylor1916}
Taylor G. I. (1916). Skin friction of the wind on the earth's surface. Proc. R. Soc. Lond. A, 92(637), 196-199.

\bibitem{usmani1994}
Usmani R. A. (1994). Inversion of a tridiagonal jacobi matrix. Linear Algebra and Its Applications, 212-213, 413-414.

\bibitem{vreugdenhil1994}
Vreugdenhil C. B. (1994). Numerical methods for shallow-water flow (Vol. 13). Springer Science \& Business Media.

\bibitem{zhang2004}
Zhang Y., Baptista A. M., \& Myers E. P. (2004). A cross-scale model for 3D baroclinic circulation in estuary–plume–shelf systems: I. Formulation and skill assessment. Continental Shelf Research, 24(18), 2187-2214. \doi{https://doi.org/10.1016/j.csr.2004.07.021}.

\end{thebibliography}

\end{document}